\documentclass[reqno]{amsart}
\usepackage{amsmath}
\usepackage{amssymb}
\usepackage{amscd}
\usepackage{graphics}
\usepackage{latexsym}
\usepackage{hyperref}

\theoremstyle{plain}
\newtheorem{theorem}{Theorem}[section]

\newtheorem{cor}[theorem]{Corollary}
\newtheorem{lem}[theorem]{Lemma}
\newtheorem{prop}[theorem]{Proposition}

\theoremstyle{definition}
\newtheorem{defn}[theorem]{Definition}

\newtheorem{rmk}[theorem]{Remark}
\newtheorem{constr}[theorem]{Construction}
\newtheorem{notat}[theorem]{Notation}
\newtheorem{ex}[theorem]{Example}
\newtheorem{hyp}[theorem]{Hypothesis}

\theoremstyle{remark}

\newcommand{\QQ}{\mathbb{Q}}

\newcommand{\ZZ}{\mathbb{Z}}

\newcommand{\AAA}{\mathbb{A}}
\newcommand{\PP}{\mathbb{P}}

\newcommand{\SP}{\text{Spec }}

\newcommand{\Kbm}[2]{\overline{\mathcal{M}}_{#1}(#2)}

\newcommand{\marpar}[1]{}

\newcommand{\lt}{\left}
\newcommand{\rt}{\right}

\newsavebox{\sembox}
\newlength{\semwidth}
\newlength{\boxwidth}

\newsavebox{\semrbox}
\newlength{\semrwidth}
\newlength{\boxrwidth}

\newcommand{\mc}{\mathcal}
\newcommand{\mf}{\mathfrak}
\newcommand{\mb}{\mathbf}

\newcommand{\OO}{\mathcal O}

\def\PP{{\mathbb P}}

\begin{document}

\title[Virtual Canonical Bundle]
{Divisor classes and the virtual canonical bundle for genus 0 maps}

\author[de Jong]{A. J. de Jong}
\address{Department of Mathematics \\
  Massachusetts Institute of Technology \\ Cambridge MA 02139}
\email{dejong@math.mit.edu}

\author[Starr]{Jason Starr}
\address{Department of Mathematics \\
  Massachusetts Institute of Technology \\ Cambridge MA 02139}
\email{jstarr@math.mit.edu}

\date{\today}

\begin{abstract}
Some divisor class relations for genus 0 curves are proved and used to
compute the Cartier divisor class of the virtual canonical bundle for
genus 0 maps to a smooth target.  Many results here first
appeared in ~\cite{QDiv} and ~\cite{Pand97}; our proofs use a
completely different method.
\end{abstract}

\maketitle

\section{Statement of results}

\noindent
Let $X$ be a smooth, projective variety over a characteristic $0$
field $k$, and denote by $\Kbm{0,r}{X,\beta}$ the Kontsevich moduli
space of genus $0$, $r$-pointed stable maps to $X$ of class $\beta$.
Behrend and Fantechi defined a \emph{perfect obstruction theory} for
$\Kbm{0,r}{X,\beta}$, i.e., a complex $E^\bullet$ perfect of amplitude
$[-1,0]$ together with a map to the cotangent complex $\phi:E^\bullet
\rightarrow L^\bullet_{\Kbm{0,r}{X,\beta}}$ such that $h^0(\phi)$ is
an isomorphism and $h^{-1}(\phi)$ is surjective.  In many cases $\phi$
is a quasi-isomorphism, and then the dualizing sheaf on
$\Kbm{0,r}{X,\beta}$ is the determinant $\text{det}(E^\bullet)$.  For
this reason $\text{det}(E^\bullet)$ is called the \emph{virtual
canonical bundle}.  This paper gives a formula,
Proposition~\ref{prop-vc}, for the virtual canonical bundle in terms
of tautological divisor classes on $X$, modulo torsion.

\medskip\noindent
Grothendieck-Riemann-Roch applies in a straightforward manner, but the
resulting formula is not optimal: it is difficult to determine whether
the virtual canonical bundle is NEF, ample, etc.  The main work in
Section~\ref{sec-div} proves divisor class relations yielding a
succinct formula for the virtual canonical bundle.  The proof reduces
to local computations for the universal family over the Artin stack of
all prestable curves of genus $0$, cf. Section~\ref{sec-local}.
Because of this, most results are stated for Artin stacks.  This leads
to one \emph{ad hoc} consruction: since there is as yet no theory of
cycle class groups for Artin stacks admitting Chern classes for all
perfect complexes of bounded amplitude, a Riemann-Roch theorem for all
perfect morphisms relatively representable by proper algebraic spaces,
and arbitrary pullbacks for all cycles coming from Chern classes, a
stand-in $Q_\pi$ is used, cf. Section~\ref{sec-Qpi} (also by avoiding
Riemann-Roch, this allows some relations to be proved ``integrally''
rather than ``modulo torsion'').  Also, although the relative Picard
functor of the universal family of genus $0$ curves is well-known, a
description is included in Sections ~\ref{sec-dec} and ~\ref{sec-cpct}
for completeness.

\medskip\noindent  
In the special case $X=\PP^n_k$, Pandharipande proved most of
these divisor class relations, ~\cite{QDiv}, and the formula for the
virtual canonical bundle, ~\cite{Pand97}, modulo numerical
equivalence.  This was certainly our inspiration, but our proofs are
completely different, yield a more general virtual canonical bundle
formula, and hold modulo torsion (and sometimes ``integrally'') rather
than modulo numerical equivalence.

\section{Decorated prestable curves} \label{sec-dec}
\marpar{sec-dec}
 
\noindent
There exists an Artin stack parametrizing prestable curves whose dual
graph is a given modular graph.  This section describes a variant
Artin stack obtained by ``decorating'' the modular graph.  Although
the variation is simple, it arises often enough to warrant a few
words.  This variant is used in the next section to describe the
closure of the identity section in the relative Picard functor for the
universal family of prestable curves of compact type.
The reference for this section is ~\cite{BM}.

\medskip\noindent
A \emph{modular graph} is a (not necessarily connected) graph $\sigma$
-- edges are undirected and \emph{tails} or half-edges are allowed --
together with a \emph{genus function} $g:\text{Vertex}(\sigma)
\rightarrow \ZZ_{\geq 0}$.  There are 2 collections of morphisms
between modular graphs: \emph{contractions} are surjective on
vertices, roughly contracting subgraphs of the domain to vertices of
the target, and \emph{graph inclusions} are inclusions of
subgraphs. (In ~\cite{BM}, the \emph{combinatorial morphisms} are
obtained from graph inclusions by adjoining formal inverses of certain
``stabilizing'' contractions; stability is not an issue here, so graph
inclusions are more appropriate).  Also, for every diagram,
$$
\begin{CD}
& & \sigma_3 \\ 
& & @VV a V \\
\sigma_2 @> \phi >> \sigma_1
\end{CD}
$$
of a contraction $\phi$ and a graph inclusion $a$, there is a
\emph{pullback diagram},
$$
\begin{CD}
\sigma_4 @> a^*\phi >> \sigma_3 \\
@V \phi^* a VV @VV a V \\
\sigma_2 @> \phi >> \sigma_1
\end{CD}
$$
of a contraction $a^*\phi$ and a graph inclusion $\phi^*a$ such that
the maps on vertices, $a\circ a^*\phi, \phi \circ
\phi^*a:\text{Vertex}(\sigma_4) \rightarrow \text{Vertex}(\sigma_1)$
are equal.  The diagram is unique up to a unique isomorphism (both as
a contraction and a graph inclusion) of $\sigma_4$.  The category of
modular graphs is denoted $\mf{G}$.

\medskip\noindent
To each prestable curve there is an associated modular graph, and to
each modular graph $\sigma$ there is an Artin stack $\mf{M}(\sigma)$
parametrizing prestable curves along with a contraction of the
associated modular graph to $\sigma$.  This defines a lax 2-functor
from $\mf{G}$ to the 2-category of Artin stacks, covariant for
contractions, contravariant for graph inclusions, and such that for
every pullback diagram there is a 2-equivalence $\mf{M}(a)\circ
\mf{M}(\phi) \Rightarrow \mf{M}(a^*\phi) \circ \mf{M}(\phi^*a)$.

\begin{defn} \label{defn-decg}
\marpar{defn-decg}
A \emph{category of decorated modular graphs} is a category $\mf{H}$
with 2 sets of morphisms -- $\mf{H}$-contractions and $\mf{H}$-graph
inclusions -- together with a functor $p:\mf{H} \rightarrow \mf{G}$
compatible with both contractions and graph inclusions satisfying the
following axioms,
\begin{enumerate}
\item[(i)]
for every $\mf{H}$-contraction $\phi:\tau_2 \rightarrow \tau_1$ and
$\mf{H}$-graph inclusion $a:\tau_3 \rightarrow \tau_1$, there exists
an object $\tau_4$, an $\mf{H}$-contraction $a^*\phi:\tau_4
\rightarrow \tau_3$ and an $\mf{H}$-graph inclusion $\phi^*a:\tau_4
\rightarrow \tau_2$ mapping under $p$ to a pullback diagram, moreover
this is unique up to unique isomorphism of $\tau_4$, and
\item[(ii)]
for every object $\tau$ in $\mf{H}$ and every contraction $\phi:p(\tau)
\rightarrow \sigma$ in $\mf{G}$, there is a $\mf{H}$-contraction $\psi$
such that $p(\psi)=\phi$, and $\psi$ is unique up to unique isomormphism.
\end{enumerate}
\end{defn}

\begin{constr} \label{constr-semigp}
\marpar{constr-semigp}
Let $A$ be an Abelian semigroup.  Define $\mf{G}_A$ to be the category
whose objects are pairs $(\sigma,\alpha)$ of a modular graph $\sigma$
together with a function $\alpha:\text{Vertex}(\sigma) \rightarrow A$,
where $\mf{G}_A$-contractions, $\phi:(\sigma_1,\alpha_1) \rightarrow
(\sigma_2,\alpha_2)$, are contractions $\phi:\sigma_1 \rightarrow
\sigma_2$ such that $\alpha_2(v) = \sum_{w\in \phi^{-1}(v)}
\alpha_1(w)$ for every $v\in \text{Vertex}(\sigma_2)$, and where
$\mf{G}_A$-graph inclusions, $a:(\sigma_1,\alpha_1) \rightarrow
(\sigma_2,\alpha_2)$, are graph inclusions $a:\sigma_1 \rightarrow
\sigma_2$ such that $\alpha_2(a(v)) = \alpha_1(v)$ for every vertex
$v\in \text{Vertex}(\sigma_1)$.  Define $p:\mf{G}_A \rightarrow
\mf{G}$ to be the obvious forgetful functor.  This is a category of
decorated modular graphs; the only one used in the rest of this paper.
\end{constr}

\medskip\noindent
The aim of this section is to construct for every object $\tau$ of
$\mf{H}$ an Artin stack $\mf{M}_{\mf{H}}(\tau)$ parametrizing
prestable curves along with a lifting of the associated modular graph
to an object of $\mf{H}$ contracting to $\tau$.  The association $\tau
\mapsto \mf{M}_{\mf{H}}(\tau)$ should define a lax 2-functor from
$\mf{H}$ to the 2-category of Artin stacks, covariant for
$\mf{H}$-contractions, contravariant for $\mf{H}$-graph inclusions,
and such that for every pullback diagram there is an associated
2-equivalence.

\begin{defn} \label{defn-sat}
\marpar{defn-sat}
Let $\mf{H}$ be a category of decorated modular graphs, considered as
a usual category whose morphisms are $\mf{H}$-contractions.  A
subcategory $\mf{H}'$ is \emph{saturated} if $\mf{H}'$ contains every
$\mf{H}$-contraction whose domain is in $\mf{H}'$.  A subcategory
$\mf{H}'$ is \emph{$p$-embedding} if the functor of categories with
contractions as morphisms, $p:\mf{H}' \rightarrow \mf{G}$, is an
equivalence to a (necessarily full) subcategory of $\mf{G}$.

\noindent
Let $\tau$ be an object of $\mf{H}$ and denote by $\mf{H}_\tau$ the
category whose objects are contractions $\phi:\tau_\phi \rightarrow
\tau$ and whose morphisms are commutative diagrams of contractions.  A
subcategory $\mf{H}'$ of $\mf{H}_\tau$ is \emph{saturated} if
$\mf{H}'$ contains every morphism in $\mf{H}_\tau$ whose domain is in
$\mf{H}'$.  A subcategory $\mf{H}'$ of $\mf{H}_\tau$ is
\emph{$p$-embedding} if the functor $p:\mf{H}' \rightarrow
\mf{G}_{p(\tau)}$ is an equivalence to a (necessarily full)
subcategory of $\mf{G}_{p(\tau)}$.  Denote by
$\text{Sat}(\mf{H}_\tau)$ the set of saturated, $p$-embedding
subcategories of $\mf{H}_\tau$ directed by reverse inclusion of
subcategories.
\end{defn}

\medskip\noindent
Let $\tau$ be an object of $\mf{H}$ and let $\mf{H}'$ be a saturated,
$p$-embedding subcategory of $\mf{H}_\tau$.  Define
$U_{\mf{H}'}(\tau)$ to be the open substack of $\mf{M}(p(\tau))$ whose
complement is the union of the images of all 1-morphisms
$\mf{M}(\phi):\mf{M}(\sigma) \rightarrow \mf{M}(p(\tau))$ such that
$\phi$ is not in the image of $p:\mf{H}' \rightarrow
\mf{G}_{p(\tau)}$.  It is straightforward that $u_{\mf{H}'}(\tau)$ is
open: the intersection with every quasi-compact open substack of
$\mf{M}(p(\tau))$ is open, and $U_{\mf{H}'}(\tau)$ is the union of
these open sets.

\medskip\noindent
Let $\mf{H}' \subset \mf{H}''$ be saturated, $p$-embedding
subcategories of $\mf{H}_\tau$.  Then $U_{\mf{H}'}(\tau) \subset
U_{\mf{H}''}(\tau)$ as subsets of $\mf{M}(p(\tau))$.  Therefore
$\mf{H}' \mapsto U_{\mf{H}'}(\tau)$ is a directed system of open
immersions of Artin stacks indexed by $\text{Sat}(\mf{H}_\tau)$.
Because this is a directed system of open immersions, the direct limit
is an Artin stack.

\begin{notat} \label{notat-MHtau}
\marpar{notat-MHtau}
Denote by $\mf{M}_{\mf{H}}(\tau)$ the direct limit of the directed
system $\mf{H}' \mapsto U_{\mf{H}'}(\tau)$.  Denote by
$\mf{M}_p(\tau): \mf{M}_{\mf{H}}(\tau) \rightarrow \mf{M}(p(\tau))$
the natural 1-morphism.  If $\mf{H}=\mf{G}_A$, also denote
$\mf{M}_{\mf{H}}(\tau)$ by $\mf{M}_A(\tau)$.
\end{notat}

\medskip\noindent
The ``points'' of $\mf{M}_{\mf{H}}(\tau)$ have a simple description.

\begin{defn} \label{defn-strict}
\marpar{defn-strict}
For every modular graph $\sigma$, define
$\mf{M}^{\text{strict}}(\sigma)$ to be the open substack of
$\mf{M}(\sigma)$ that is the complement of the images of all
$\mf{M}(\phi)$ where $\phi:\sigma' \rightarrow \sigma$ is a
non-invertible contraction.
\end{defn}

\begin{lem} \label{lem-pts}
\marpar{lem-pts}
Let $\tau$ be an object of $\mf{H}$ and let $\phi:\sigma \rightarrow
p(\tau)$ be a contraction.  The 2-fibered product
$\mf{M}^{\text{strict}}(\sigma) \times_{\mf{M}(p(\tau))}
\mf{M}_{\mf{H}}(\tau)$ is equivalent to a disjoint union of copies of
$\mf{M}^{\text{strict}}(\sigma)$ indexed by equivalence classes of
contractions $\psi$ in $\mf{H}_\tau$ such that $p(\psi)=\phi$.
\end{lem}

\begin{proof}
Let $(\eta,\zeta,\theta)$ be an object of the 2-fibered product, i.e.,
a triple of an object of $\mf{M}^{\text{strict}}(\sigma)$, an object
of $\mf{M}_{\mf{H}}(\tau)$ and an equivalence
$\theta:\mf{M}(\phi)(\eta) \rightarrow \mf{M}_{p}(\tau)(\zeta)$.
There is a saturated, $p$-embedding subcategory $\mf{H}'$ such that
$\mf{M}_p(\tau)(\zeta)$ is in $U_{\mf{H}'}(\tau)$.  Because this is in
the image of $\mf{M}(\phi)$, there is a contraction $\psi:\tau_\psi
\rightarrow \tau$ in $\mf{H}'$ such that $\phi=p(\psi)$.  Because
$\mf{H}'$ is $p$-embedding, $\psi$ is unique up to unique isomorphism.
By the nature of the direct limit, $\psi$ is independent of the choice
of $\mf{H}'$.

\medskip\noindent
Conversely, given an object $\eta$ of $\mf{M}^{\text{strict}}(\sigma)$
and a contraction $\psi:\tau_\psi \rightarrow \tau$ in $\mf{H}_\tau$
such that $p(\psi) = \phi$, define $\mf{H}'$ to be the subcategory of
$\mf{H}_\tau$ consisting of all contractions through which $\psi$
factors.  By Definition~\ref{defn-decg}(ii), this is a saturated,
$p$-embedding subcategory.  And $\mf{M}(\phi)(\eta)$ is in
$U_{\mf{H}}(\tau)$.  The image in the direct limit is an object
$\zeta$, and there is a canonical isomorphism $\theta:
\mf{M}(\phi)(\eta) \rightarrow \mf{M}_p(\tau)(\zeta)$.  Thus
$(\eta,\zeta,\theta)$ is an object of the 2-fibered product.

\medskip\noindent
It is left to the reader to verify these operations give an
equivalence of stacks.
\end{proof}

\medskip\noindent
\textbf{Note:}  
The functorialities are only sketched.  Given an $\mf{H}$-contraction
$\phi:\tau_1 \rightarrow \tau_2$, Definition~\ref{defn-decg}(ii) gives
a map of directed sets $\text{Sat}(\phi):\text{Sat}(\mf{H}_{\tau_1})
\rightarrow \text{Sat}(\mf{H}_{\tau_1})$, and composition with
$\mf{M}(p(\phi))$ gives a compatible family of 1-morphisms of directed
systems.  This defines the 1-morphism $\mf{M}_{\mf{H}}(\phi)$.  Given
an $\mf{H}$-graph inclusion $a:\tau_1 \rightarrow \tau_2$, existence
of pullback diagrams, Definition~\ref{defn-decg}(i), gives a map of
directed sets $\text{Sat}(a):\text{Sat}(\mf{H}_{\tau_2}) \rightarrow
\text{Sat}(\mf{H}_{\tau_1})$, and composition with $\mf{M}(p(a))$
gives a compatible family of 1-morphisms of directed systems.  This
defines the 1-morphism $\mf{M}_{\mf{H}}(a)$.  The rest of the
compatibilities are straightforward.

\section{The universal relative Picard for curves of compact type}
\label{sec-cpct}
\marpar{sec-cpct}

\noindent
The results in this section are well-known, and easily follow from
~\cite{Raynaud} and ~\cite{Ner}.  It is useful in the rest of the
paper to gather the results here.

\begin{notat} \label{notat-H}
\marpar{notat-H}
Denote by $\mf{H} \subset \mf{G}_\ZZ$ the full subcategory of objects
$(\sigma,\alpha)$ such that $\sigma$ is a forest of trees, i.e., the
graph has no cycles.  For each triple of integers $g,n \geq 0$ and
$e$, denote by $\tau_{g,n}(e)$ the object of $\mf{H}$ consisting of a
tree $\sigma_{g,n}$ with a single vertex of genus $g$ and $n$ flags,
such that $\alpha(v)=n$.
\end{notat}

\medskip\noindent
Denote by $\pi:\mc{C} \rightarrow \mf{M}_{\mf{H}}(\tau_{g,0}(0))$ the
pullback from $\mf{M}(\sigma_{g,0})$ of the universal curve.  For each
4-tuple $A=((g',g''),(e',e''))$ of integers $g',g''\geq 0$,
$g'+g''=g$, and integers $e',e''$, $e'+e''=0$, denote by $\tau_A$ the
tree with vertices $v',v''$ such that $g(v')=g', \alpha(v')=e'$ and
$g(v'')=g'', \alpha(v'')=e''$.  Denote by $\phi:\tau_A \rightarrow
\tau_{g,0}(0)$ the canonical contraction.  The 2-fibered product
$\mf{M}_{\mf{H}}(\tau_A) \times_{\mf{M}_{\mf{H}}(\tau_{g,0}(0))}
\mc{C}$ has 2 irreducible components $\mc{C}', \mc{C}''$ corresponding
to the vertices $v', v''$.  There is a unique effective Cartier
divisor $\mc{D} \subset \mc{C}$ such that for every
$A=((g',g''),(e',e''))$,
$$
\mf{M}^{\text{strict}}_{\mf{H}}(\tau_A)
\times_{\mf{M}_{\mf{H}}(\tau_{g,0}(0))} \mc{D}
$$ 
is empty if $e'=e''=0$ and is $e'\mc{C}''$ if
$e'>0$.  

\medskip\noindent
Let $U\subset \mf{M}(\sigma_{g,0})$ denote the open substack that is
the image of $\mf{M}_p(\tau_{g,0}(0))$, i.e., $U$ is the Artin stack
of $n$-pointed, genus $g$ curves of \emph{compact type}.  The
1-morphism $\pi:\mc{C}_U \rightarrow U$ is cohomologically flat, so by
~\cite[Thm. 7.3]{Artin} the relative Picard functor of the universal
curve over $U$ is a 1-morphism $\text{pr}:\text{Pic}_{\mc{C}_U/U}
\rightarrow U$ relatively representable by \emph{non-separated}
algebraic spaces.  The closure $E_{\mc{C}_U/U}$ of the identity
section gives a closed substack of $\text{Pic}_{\mc{C}_U/U}$ which is
relatively representable over $U$ by non-separated group algebraic
spaces, ~\cite[Prop. 5.2]{Raynaud}.  The quotient $Q_{\mc{C}_U/U}$ of
$\text{Pic}_{\mc{C}_U/U}$ by $E_{\mc{C}_U/U}$ is a stack that is
relatively representable over $U$ by a countable disjoint union of
smooth, proper group algebraic spaces, ~\cite[Thm. 4.1.1]{Raynaud}
(properness requires a bit more, see ~\cite[Ex. 8, p. 246]{Ner}).  The
next lemma describes $E_{\mc{C}_U/U}$.

\medskip\noindent
The invertible sheaf $\OO_{\mc{C}}(\mc{D})$ defines a
1-morphism $f:\mf{M}_{\mf{H}}(\tau_{g,0}(0)) \rightarrow
\text{Pic}_{\mc{C}_U/U}$, and there is a natural 2-equivalence of
$\text{pr}\circ f$ with $\mf{M}_{\mf{H}}(\tau_{g,0}(0))$.

\begin{lem} \label{lem-subgp} 
\marpar{lem-subgp}
The 1-morphism $f$ defines an equivalence to $E_{\mc{C}_U/U}$, the
closure of the 
identity section of $\text{Pic}_{\mc{C}_U/U}$.  
Denoting by
$Q^0_{\mc{C}_U/U}$ the identity component of the quotient and by
$\text{Pic}^0_{\mc{C}_U/U}$ the preimage, there are 1-morphisms
$$
\text{Pic}^0_{\mc{C}_U/U} \rightleftarrows
\mf{M}_{\mf{H}}(\tau_{g,0}) \times
Q^0_{\mc{C}_U/U}
$$ 
giving an equivalence of stacks
over $U$, and splitting the extension of group algebraic spaces over
$U$.
\end{lem}  

\begin{proof}
It is easy to see $f$ is an equivalence to its image which is a
subgroup of $E_{\mc{C}_U/U}$.  To prove the image of $f$ is all of
$E_{\mc{C}_U/U}$, by the valuative criterion of closedness it suffices
to check equality of pullbacks for every map of a DVR to
$\mf{M}_{g,0}$ sending the generic point to
$\mf{M}_{g,0}^{\text{strict}}$. By ~\cite[Prop. 6.1.3]{Raynaud}, the
sections of $E_{\mc{C}_U/U}$ over are a DVR are just the quotient of
the free Abelian group on the irreducible components of the closed
fiber by the subgroup generated by the entire fiber.  By
Lemma~\ref{lem-pts} the same is true for the pullback of
$\mf{M}_{\mf{H}}(\tau_{g,0}(0))$, and it is clear the map between them
is an isomorphism.

\medskip\noindent
The splitting of $\text{Pic}^0_{\mc{C}_U/U} \rightarrow
Q^0_{\mc{C}_U/U}$ is given by the subfunctor of
$\text{Pic}^0_{\mc{C}_U/U}$ of invertible sheaves whose degree on
every irreducible component of every fiber is $0$, denoted by $P^0$ in
\cite{Raynaud}.
\end{proof}

\medskip\noindent
In the special case that $g=0$, more is true.  First of all,
$\mf{M}_{\ZZ}(\tau) = \mf{M}_{\mf{H}}(\tau)$ for every $\tau$ of genus
$0$.  Secondly,  $U=\mf{M}_{0,n}$.  

\begin{cor}[Raynaud, Prop. 9.3.1, ~\cite{Raynaud}] \label{cor-subgp}
\marpar{cor-subgp}
For $g=0$, $\text{Pic}^0_{\mc{C}/\mf{M}_{0,0}}$ is equivalent
to $\mf{M}_{\ZZ}(\tau_{0,0}(0))$.
\end{cor}

\medskip\noindent
Moreover, the union $\cup_{e\in \ZZ} \mf{M}_{\ZZ}(\tau_{0,0}(e))$ is a
group algebraic space over $\mf{M}_{0,0}$ containing
$\mf{M}_{\ZZ}(\tau_{0,0}(0))$ as a subgroup algebraic space over
$\mf{M}_{0,0}$.  Essentially, given a contraction
$\phi:\sigma\rightarrow \sigma_{0,0}$ and given liftings
$\psi_i:(\sigma,\alpha_i) \rightarrow \tau_{0,0}(e_i)$ for $i=1,2$,
addition is determined by $\psi_1+\psi_2=\psi:(\sigma,
\alpha_1+\alpha_2) \rightarrow \tau_{0,0}(e_1+e_2)$.  The \emph{total
degree map} gives an isomorphism of $\text{Pic}/\text{Pic}^0$ with
$\ZZ$.  The following result is easy.

\begin{lem} \label{lem-g0}
\marpar{lem-g0}
For $g=0$, for every $e$ there is an equivalence of stacks over
$\mf{M}_{0,0}$, $\mf{M}_{\ZZ}(\tau_{0,0}(e)) \rightleftarrows
\text{Pic}^e_{\mc{C}/\mf{M}_{0,0}}$, such that the equivalence
$\cup_{e\in \ZZ} \mf{M}_{\ZZ}(\tau_{0,0}(e)) \rightleftarrows
\text{Pic}_{\mc{C}/\mf{M}_{0,0}}$ is an equivalence of group algebraic
spaces over $\mf{M}_{0,0}$ and is compatible with the equivalence in
Corollary~\ref{cor-subgp}.
\end{lem}

\subsection{Notation for boundary divisor classes}
\label{subsec-bound}
\marpar{subsec-bound}

\noindent
Let $r\geq 0$ be an integer, and let $e_1,\dots,e_r$ be integers.
Denote by $\mf{M}_{\ZZ}(\tau_{0,0}(e_1,\dots,e_r))$ the 2-fibered
product,
$$
\mf{M}_{\ZZ}(\tau_{0,0}(e_1)) \times_{\mf{M}_{0,0}}
\mf{M}_{\ZZ}(\tau_{0,0}(e_2)) \times_{\mf{M}_{0,0}} \dots
\times_{\mf{M}_{0,0}} \mf{M}_{\ZZ}(\tau_{0,0}(e_r)).
$$
For each $i=1,\dots,r$, let $(e'_i,e''_i)$ be a pair of integers such
that $e'_i + e''_i = e_i$.  Let $\sigma$ be the modular graph with two
vertices $v',v''$ with $g(v')=g(v'')=0$, one edge connecting $v',
v''$, and no tails.  Let $\phi:\sigma \rightarrow \sigma_{0,0}$ be the
canonical contraction.  Denote by
$\zeta:\mf{M}^{\text{strict}}(\sigma) \rightarrow
\mf{M}_{\ZZ}(\tau_{g,0}(e_1,\dots,e_r))$ the 1-morphism whose
projection to the $i^{\text{th}}$ factor is determined via
Lemma~\ref{lem-pts} by the lifting $\psi_i:(\sigma,\alpha_i)
\rightarrow \tau_{g,0}(e_i)$ of $\phi$ such that $\alpha_i(v') = e'_i,
\alpha_i(v'')=e_i''$.  Define $\Delta_{(e'_1,e''_1,\dots,e'_r,e''_r)}$
to be the effective Cartier divisor on
$\mf{M}_{\ZZ}(\tau_{0,0}(e_1,\dots,e_r))$ that is the closure of the
image of $\zeta$.

\medskip\noindent
Let $\pi:C\rightarrow M$ be a flat 1-morphism relatively represented
by proper algebraic spaces whose geometric fibers are connected,
at-worst-nodal curves of arithmetic genus $0$.  Let $D_1,\dots,D_r$ be
Cartier divisor classes on $C$ of relative degrees $e_1,\dots,e_r$.
Let $f(e'_1,e''_1,\dots,e'_r,e''_r)$ be a function on $\ZZ^{2r}$ with
values in $\ZZ$, resp. $\QQ$, etc.  Denote by $\xi: M \rightarrow
\mf{M}_{\ZZ}(\tau_{0,0}(e_1,\dots,e_r))$ the 1-morphism whose
projection to the $i^\text{th}$ factor,
$\text{Pic}^{e_i}_{\mf{C}/\mf{M}_{0,0}}$ is determined by
$\OO_C(D_i)$.

\begin{notat} \label{notat-boundary}
\marpar{notat-boundary}
Denote by,
$$
\sum_{(\beta',\beta'')} f(\langle D_1,\beta' \rangle, \langle
D_1, \beta'' \rangle, \dots, \langle D_r, \beta' \rangle, \langle D_r,
\beta'' \rangle) \Delta_{\beta',\beta''}
$$
the Cartier divisor class, resp. $\QQ$-Cartier divisor class, etc., 
that is the pullback by $\xi$ of the Cartier divisor class, etc.,
$$
\sum_{(e'_1,e''_1,\dots,e'_r,e''_r)} f(e'_1,e''_1,\dots,e'_r,e''_r)
\Delta_{(e'_1,e''_1,\dots,e'_r,e''_r)},
$$
the summation over all sequences $(e'_1,e''_1,\dots,e'_r,e''_r)$ with
$e'_i + e''_i = e_i$.  
If $f(e'_1,e''_1,\dots,e'_r,e''_r) = f(e''_1,e'_1,\dots,e''_r,e'_r)$,
denote by,
$$
{\sum_{(\beta',\beta'')}}^\prime f(\langle D_1,\beta' \rangle, \langle
D_1, \beta'' \rangle, \dots, \langle D_r, \beta' \rangle, \langle D_r,
\beta'' \rangle) \Delta_{\beta',\beta''}
$$
the pullback by $\xi$ of,
$$
{\sum_{(e'_1,e''_1,\dots,e'_r,e''_r)}}^\prime
f(e'_1,e''_1,\dots,e'_r,e''_r)
\Delta_{(e'_1,e''_1,\dots,e'_r,e''_r)},
$$
where the summation is over equivalence classes of sequences
$(e'_1,e''_1,\dots,e'_r,e''_r)$ such that $e'_i+e''_i=e_i$ 
under the equivalence relation
$(e'_1,e''_1,\dots,e'_r,e''_r) \sim (e''_1,e'_1,\dots,e''_r,e'_r)$.  
\end{notat}

\begin{ex} \label{ex-boundary}
\marpar{ex-boundary}
Let $n\geq 0$ be an integer and let $(A,B)$ be a partition of
$\{1,\dots,n\}$.  For the universal family over $\mf{M}_{0,n}$, denote
by $s_1,\dots,s_n$ the universal sections.  Then,
$$
\sum_{\beta',\beta''} \prod_{i\in A} \langle s_i,\beta' \rangle
\cdot \prod_{j\in B} \langle s_j, \beta'' \rangle
\Delta_{\beta',\beta''} 
$$
is the Cartier divisor class of the boundary divisor $\Delta_{(A,B)}$.  
\end{ex}

\section{The functor $Q_\pi$} \label{sec-Qpi}
\marpar{sec-Qpi}

\noindent
Let $M$ be an Artin stack, and let $\pi:C \rightarrow M$ be a flat
1-morphism, relatively representable by proper algebraic spaces whose
geometric fibers are connected, at-worst-nodal curves of arithmetic
genus $0$.  There exists an invertible dualizing sheaf $\omega_\pi$,
and the relative trace map, $\text{Tr}_\pi: R\pi_* \omega_\pi[1]
\rightarrow \OO_M$ is a quasi-isomorphism.  In particular,
$\text{Ext}^1_{\OO_C}(\omega_\pi,\OO_C)$ is canonically isomorphic to
$H^0(M,\OO_M)$.  Therefore $1\in H^0(M,\OO_M)$ determines an extension
class, i.e., a short exact sequence,
$$
\begin{CD}
0 @>>> \omega_\pi @>>> E_\pi @>>> \OO_C @>>> 0.
\end{CD}
$$
The morphism $\pi$ is perfect, so for every complex $F^\bullet$ perfect
of bounded amplitude on $C$, $R\pi_* F^\bullet$ is a perfect complex of
bounded amplitude on $M$.  By ~\cite{detdiv}, the
determinant of a perfect complex of bounded amplitude is defined.

\begin{defn} \label{defn-E}
\marpar{defn-E}
For every complex $F^\bullet$ perfect of bounded amplitude on $C$, 
define $Q_\pi(F^\bullet) = \text{det}(R\pi_* E_\pi\otimes F^\bullet)$.
\end{defn}

\medskip\noindent
There is another interpretation of $Q_\pi(F^\bullet)$.

\begin{lem} \label{lem-interp}
\marpar{lem-interp}
For every complex $F^\bullet$ perfect of bounded amplitude on $C$, 
$$
Q_\pi(F^\bullet) \cong
\text{det}(R\pi_*(F^\bullet)) \otimes 
\text{det}(R\pi_*((F^\bullet)^\vee))^\vee.
$$
\end{lem}

\begin{proof}
By the short exact sequence for $E_\pi$, $Q_\pi(F^\bullet) \cong
\text{det}(R\pi_*(F^\bullet)) \otimes \text{det}(R\pi_*(\omega_\pi
\otimes F^\bullet))$.  The lemma follows by duality.
\end{proof}

\medskip\noindent
It is straightforward to compute $F^\bullet$ whenever there exist
cycle class groups for $C$ and $M$ such that Chern classes are defined
for all perfect complexes of bounded amplitude and such that
Grothendieck-Riemann-Roch holds for $\pi$.

\begin{lem} \label{lem-GRR}
\marpar{lem-GRR}
If there exist cycle class groups for $C$ and $M$ such that Chern
classes exist for all perfect complexes of bounded amplitude and such
that Grothendieck-Riemann-Roch holds for $\pi$, then modulo $2$-power
torsion, the first Chern class of $Q_\pi(F^\bullet)$ is
$\pi_*(C_1(F^\bullet)^2 - 2C_2(F^\bullet))$.
\end{lem}

\begin{proof}
Denote the Todd class of $\pi$ by $\tau = 1 + \tau_1 + \tau_2 +
\dots$.  Of course $\tau_1 = -C_1(\omega_\pi)$.  By GRR,
$\text{ch}(R\pi_* \OO_C) = \pi_*(\tau)$.  The canonical map $\OO_M
\rightarrow R\pi_*\OO_C$ is a quasi-isomorphism.  Therefore
$\pi_*(\tau_2)=0$, modulo $2$-power torsion.  By additivity of the
Chern character, $\text{ch}(E_\pi) = 2 + C_1(\omega_\pi) +
\frac{1}{2}C_1(\omega_\pi)^2 + \dots$.  Therefore,
$$
\text{ch}(E_\pi)\cdot \tau = 2 + 2\tau_2 + \dots
$$
So for any complex $F^\bullet$ perfect of bounded amplitude,
$$
\begin{array}{c}
\text{ch}(E_\pi\otimes F^\bullet)\cdot \tau = \text{ch}(F^\bullet) \cdot
\text{ch}(E_\pi) \cdot \tau = \\
(\text{rk}(F^\bullet) + C_1(F^\bullet) +
\frac{1}{2}(C_1(F^\bullet)^2 - 2C_2(F^\bullet)) +\dots)(2 + 2\tau_2+
\dots). 
\end{array}
$$
Applying $\pi_*$ gives,
$$
2\pi_*(C_1(F^\bullet)) + \pi_*(C_1(F^\bullet)^2 - 2C_2(F^\bullet)) + \dots
$$
Therefore the first Chern class of $\text{det}(R\pi_*(E_\pi\otimes
F^\bullet))$ is $\pi_*(C_1(F^\bullet)^2 - 2C_2(F^\bullet))$, modulo
$2$-power torsion.
\end{proof}

\begin{rmk} \label{rmk-Q}
The point is this.  In every reasonable case, $Q_\pi$ is just
$\pi_*(C_1^2-2C_2)$.  Moreover $Q_\pi$ is compatible with base-change
by arbitrary 1-morphisms.  This allows to reduce certain computations
to the Artin stack of all genus $0$ curves.  As far as we are aware,
no one has written a definition of cycle class groups for all locally
finitely presented Artin stacks that has Chern classes for all perfect
complexes of bounded amplitude, has pushforward maps and
Grothendieck-Riemann-Roch for perfect 1-morphisms representable by
proper algebraic spaces, and has pullback maps by arbitrary
1-morphisms for cycles coming from Chern classes.  Doubtless such a
theory exists; whatever it is, $Q_\pi = \pi_*(C_1^2-2C_2)$.
\end{rmk}

\medskip\noindent
Let the following diagram be 2-Cartesian,
$$
\begin{CD}
C' @> \zeta_C >> C \\
@V \pi' VV @VV \pi V \\
M' @> \zeta_M >> M
\end{CD}
$$
together with a 2-equivalence $\theta:\pi\circ \zeta_C \Rightarrow
\zeta_M \circ \pi'$.  

\begin{lem} \label{lem-pullback}
\marpar{lem-pullback}
For every complex $F^\bullet$ perfect of bounded amplitude on $C$, 
$\zeta_M^* Q_\pi(F^\bullet)$ is isomorphic to $Q_{\pi'}(\zeta_C^*
F^\bullet)$.
\end{lem}

\begin{proof}
Of course $\zeta_C^* E_\pi = E_{\pi'}$.  And $\zeta_M^* R\pi_*$ is
canonically equivalent to $R(\pi')_* \zeta_C^*$ for perfect
complexes of bounded amplitude.  Therefore $\zeta_M^*
Q_\pi(F^\bullet)$ equals $\text{det}(\zeta_M^* R\pi_*(E_\pi\otimes
F^\bullet))$ equals $\text{det}(R(\pi')_* \zeta_C^*(E_\pi \otimes
F^\bullet))$ equals $\text{det}(R(\pi')_* E_{\pi'}\otimes \zeta_C^*
F^\bullet)$ equals $Q_{\pi'}(\zeta_C^* F^\bullet)$.
\end{proof}

\begin{lem} \label{lem-inv}
\marpar{lem-inv}
Let $L$ be an invertible sheaf on $C$ of relative degree $e$ over
$M$.  For every invertible sheaf $L'$ on $M$, $Q_\pi(L\otimes \pi^*L')
\cong Q_\pi(L)\otimes (L')^{2e}$.  In particular, if $e=0$,
$Q_\pi(L\otimes \pi^* L') \cong Q_\pi(L)$.
\end{lem}

\begin{proof}
To compute the rank of $R\pi_*(E_\pi\otimes F^\bullet)$ over any
connected component of $M$, it suffices to base-change to the spectrum
of a field mapping to that component.  Then, by
Grothendieck-Riemann-Roch, the rank is $2\text{deg}(C_1(F^\bullet))$.
In particular, $R\pi_*(E_\pi \otimes L)$ has rank $2e$.  

\medskip\noindent
By the projection formula, $R\pi_*(E_\pi\otimes L\otimes \pi^* L')
\cong R\pi_*(E_\pi \otimes L)\otimes L'$.  Of course
$\text{det}(R\pi_*(E_\pi \otimes L)\otimes L') = Q_\pi(L)\otimes
(L')^\text{rank}$.  This follows from the uniqueness of $\text{det}$:
for any invertible sheaf $L'$ the association $F^\bullet \mapsto
\text{det}(F^\bullet \otimes L')\otimes
(L')^{-\text{rank}(F^\bullet)}$ also satisfies the axioms for a
determinant function and is hence canonically isomorphic to
$\text{det}(F^\bullet)$.  Therefore $Q_\pi(L\otimes \pi^*L') =
Q_\pi(L)\otimes (L')^{2e}$.
\end{proof}

\section{Local computations} \label{sec-local}
\marpar{sec-local}

\noindent
This section contains 2 computations: $Q_\pi(\omega_\pi)$ and
$Q_\pi(L)$ for every invertible sheaf on $C$ of relative degree $0$.
Because of Lemma~\ref{lem-pullback} the first computation reduces to
the universal case over $\mf{M}_{0,0}$.  Because of
Lemma~\ref{lem-pullback} and Lemma~\ref{lem-inv}, the second
compuation reduces to $\OO_{\mc{C}}(\mc{D})$ over
$\mf{M}_{\ZZ}(\tau_{0,0}(0))$.  In each case the computation is
performed locally.

\subsection{Computation of $Q_\pi(\omega_\pi)$} \label{subsec-Qomega}
\marpar{subsec-Qomega}
Associated to $\pi_C:C\rightarrow M$, there is a 1-morphism $\zeta_M:M
\rightarrow \mf{M}_{0,0}$, a 1-morphism $\zeta_C: C \rightarrow
\mc{C}$, and a 2-equivalence $\theta:\pi_{\mc{C}} \circ \zeta_C
\Rightarrow \zeta_M\circ \pi_C$ such that the following diagram is
2-Cartesian,
$$
\begin{CD}
C @> \zeta_C >> \mc{C} \\
@V \pi_C VV @VV \pi_{\mc{C}} V \\
M @> \zeta_M >> \mf{M}_{0,0}
\end{CD}
$$
Of course $\omega_{\pi_C}$ is isomorphic to $\zeta_C^*
\omega_{\pi_\mc{C}}$.  By Lemma~\ref{lem-pullback},
$Q_{\pi_C}(\omega_{\pi_C}) \cong \zeta_M^*
Q_{\pi_{\mc{C}}}(\omega_{\pi_{\mc{C}}})$.  So the computation of
$Q_{\pi_C}(\omega_{\pi_C})$ is reduced to the universal family.

\medskip\noindent
Let the open substack $U_1\subset \mf{M}_{0,0}$ be the complement of
the union of the images of $\mf{M}(\phi):\mf{M}(\sigma) \rightarrow
\mf{M}_{0,0}$ as $\phi:\sigma \rightarrow \sigma_{0,0}$ ranges over
all contractions such that $\#\text{Vertex}(\sigma) \geq 3$.  Let $U_2
\subset U_1$ be the open substack
$\mf{M}^{\text{strict}}(\sigma_{0,0})$.

\begin{prop} \label{prop-Qomega}
\marpar{prop-Qomega}
\begin{enumerate}
\item[(i)]
Over the open substack $U_1$, $\omega_\pi^\vee$ is $\pi$-relatively
ample.  
\item[(ii)]
Over $U_1$,
$R^1\pi_*\omega_\pi^\vee|_{U_1} = (0)$ and $\pi_*
\omega_\pi^\vee|_{U_1}$ is locally free of rank 3.  
\item[(iii)]
Over $U_2$, there is
a canonical isomorphism $i:\text{det}(\pi_* \omega_\pi^\vee|_{U_2})
\rightarrow \OO_{U_2}$. 
\item[(iv)]
The image of $\text{det}(\pi_*\omega_\pi^\vee|_{U_1}) \rightarrow
\text{det}(\pi_*\omega_\pi^\vee|_{U_2}) \xrightarrow{i} \OO_{U_2}$ is
$\OO_{U_1}(-\Delta) \subset \OO_{U_2}$.
\item[(v)]
Over $U_1$, $Q_\pi(\omega_\pi)|_{U_1} \cong \OO_{U_1}(-\Delta)$.
Therefore on all of $\mf{M}_{0,0}$, $Q_\pi(\omega_\pi) \cong
\OO_{\mf{M}_{0,0}}(-\Delta)$.  
\end{enumerate}
\end{prop}

\begin{proof}
Over $\ZZ$, let $V = \ZZ\{\mb{e}_0,\mb{e}_1\}$ be a free module of
rank $2$.  Choose dual coordinates $y_0,y_1$ for $V^\vee$.  Let
$\PP^1_\ZZ = \PP(V)$ be the projective space with homogeneous
coordinates $y_0, y_1$.  Let $\AAA^1_\ZZ$ be the affine space with
coordinate $x$.  Denote by $Z\subset \AAA^1_\ZZ \times \PP^1_\ZZ$ the
closed subscheme $\mathbb{V}(x,y_1)$, i.e., the image of the section
$(0,(1,0))$.  Let $\nu:C\rightarrow \AAA^1_\ZZ \times \PP^1_\ZZ$ be
the blowing up of $Z$.  Denote by $E\subset C$ the exceptional
divisor.

\medskip\noindent
Define $\pi:C\rightarrow \AAA^1_\ZZ$ to be $\text{pr}_{\AA^1}\circ
\nu$.  This is a flat, proper morphism whose geometric fibers are
connected, at-worst-nodal curves of arithmetic genus $0$.  Moreover,
no geometric fiber has more than 1 node.  Thus there is a 1-morphism
$\zeta:\AAA^1_\ZZ \rightarrow U_1$ such that the pullback of $\mc{C}$
is equivalent to $C$.  It is straightforward that $\zeta$ is smooth
and is surjective on geometric points.  Thus (i) and (ii) can be
checked after base-change by $\zeta$.  Also (iv) will reduce to a
computation after base-change by $\zeta$.

\medskip\noindent
\textbf{(i) and (ii):} Denote by $\PP^2_\ZZ$ the projective space with
coordinates $u_0,u_1,u_2$.  There is a rational transformation
$f:\AAA^1_\ZZ \times \PP^1_\ZZ \dashrightarrow \AAA^1_\ZZ \times
\PP^2_\ZZ$ by
$$
\begin{array}{ccc}
f^*x & = & x, \\
f^*u_0 & = & xy_0^2, \\
f^*u_1 & = & y_0y_1, \\
f^*u_2 & = & y_1^2
\end{array}
$$
By local computation, this extends to a morphism $f:C\rightarrow
\AAA^1_\ZZ \times \PP^2_\ZZ$ that is a closed immersion and whose
image is $\mathbb{V}(u_0u_2-xu_1^2)$.  By the adjunction formula,
$\omega_\pi$ is the pullback of $\OO_{\PP^2}(-1)$.  In particular,
$\omega_\pi^\vee$ is very ample. Moreover, because
$H^1(\PP^2_\ZZ,\OO_{\PP^2}(-1)) = (0)$, also $H^1(C,\omega_\pi^\vee) =
(0)$.  By cohomology and base-change results,
$R^1\pi_*(\omega_\pi^\vee) = (0)$ and $\pi_*(\omega_\pi^\vee)$ is
locally free of rank 3.  

\medskip\noindent
\textbf{(iii):} The curve $\PP^1_\ZZ = \PP(V)$ determines a morphism
$\eta: \SP(\ZZ) \rightarrow U_2$.  This is smooth and surjective on
geometric points.  Moreover it gives a realization of $U_2$ as the
classifying stack of the group scheme $\text{Aut}(\PP(V)) =
\textbf{PGL}(V)$.  Taking the exterior power of the Euler exact
sequence, $\omega_{\PP(V)/\ZZ} = \bigwedge^2(V^\vee)\otimes
\OO_{\PP(V)}(-2)$.  Therefore $H^0(\PP(V),\omega_{\PP(V)/\ZZ}^\vee)$
equals $\bigwedge^2(V)\otimes \text{Sym}^2(V^\vee)$ as a
representation of $\textbf{GL}(V)$.  The determinant of this
representation is the trivial character of $\textbf{GL}(V)$.
Therefore it is the trivial character of $\textbf{PGL}(V)$.  This
gives an isomorphism of $\text{det}(\pi_*\omega_\pi|_{U_2})$ with
$\OO_{U_2}$.

\medskip\noindent
\textbf{(iv):} This can be checked after pulling back by $\zeta$.  The
pullback of $U_2$ is $\mathbb{G}_{m,\ZZ} \subset \AAA^1_\ZZ$.  The
pullback of $i$ comes from the determinant of
$H^0(\mathbb{G}_{m,\ZZ}\times \PP^1_\ZZ,\omega_\pi^\vee) =
\bigwedge^2(V) \otimes \text{Sym}^2(V^\vee) \otimes
\OO_{\mathbb{G}_m}$.  By the adjunction formula, $\omega_{C/\AAA^1} =
\nu^*\omega_{\AAA^1\times \PP^1/\AAA^1}(E)$.  Hence
$\nu_*\omega_{C/\AAA^1}^\vee = I_Z \omega_{\AAA^1\times
\PP^1/\AAA^1}$.  Therefore the canonical map,
$$
H^0(C,\omega_{C/\AAA^1}^\vee) \rightarrow H^0(\AAA^1_\ZZ \times
\PP^1_\ZZ, \omega^\vee_{\AAA^1\times \PP^1/\AAA^1}),
$$ 
is given by,
$$
\begin{array}{l}
\OO_{\AAA^1}\{\mb{f}_0,\mb{f}_1,\mb{f}_2\} \rightarrow
\bigwedge^2(V)\otimes \text{Sym}^2(V^\vee) \otimes \OO_{\AAA^1}, \\
\mb{f}_0 \mapsto x\cdot (\mb{e}_0\wedge\mb{e}_1)\otimes y_0^2, \\
\mb{f}_1 \mapsto  \ \ \ \     (\mb{e}_0\wedge \mb{e}_1) \otimes y_0y_1, \\
\mb{f}_2 \mapsto  \ \ \ \      (\mb{e}_0\wedge \mb{e}_1) \otimes y_1^2
\end{array}
$$
It follows that $\text{det}(\pi_*\omega_\pi^\vee) \rightarrow
\OO_{\mathbb{G}_m}$ has image $\langle x \rangle \OO_{\AAA^1}$, i.e.,
$\eta^* \OO_{U_1}(-\Delta)$.

\medskip\noindent
\textbf{(v):} By the short exact sequence for $E_\pi$,
$Q_\pi(\omega_\pi) = \text{det}(R\pi_*\omega_\pi)\otimes
\text{deg}(R\pi_* \omega_\pi^2)$.  Because the trace map is a
quasi-isomorphism, $\text{det}(R\pi_* \omega_\pi) = \OO_{U_1}$.  By
(ii) and duality, $$ \text{det}(R\pi_* \omega_\pi^2) \cong
\text{det}(R^1\pi_*\omega_\pi^2)^vee \cong \text{det}(\pi_*
\omega_\pi^\vee).
$$  
By (iv), this is $\OO_{U_1}(-\Delta)$.  Therefore $Q_\pi(\omega_\pi)
\cong \OO_{U_1}(-\Delta)$ on $U_1$.  Because $\mf{M}_{0,0}$ is
regular, and because the complement of $U_1$ has codimension $2$, this
isomorphism of invertible sheaves extends to all of $\mf{M}_{0,0}$.
\end{proof}

\medskip\noindent
The sheaf of relative differentials $\Omega_\pi$ is a pure coherent
sheaf on $\mc{C}$ of rank $1$, flat over $\mf{M}_{0,0}$ and is
quasi-isomorphic to a perfect complex of amplitued $[-1,0]$.  

\begin{lem} \label{lem-Omega}
\marpar{lem-Omega}
The perfect complex $R\pi_*\Omega_\pi$ has rank $-1$ and determinant
$\cong \OO_{\mf{M}_{0,0}}(-\Delta)$.  The perfect complex $R\pi_*
R\textit{Hom}_{\OO_{\mc{C}}}(\Omega_\pi,\OO_{\mc{C}})$ has rank $3$
and determinant $\cong \OO_{\mf{M}_{0,0}}(-2\Delta)$.
\end{lem}

\begin{proof}
There is a canonical injective sheaf homomorphism $\Omega_\pi
\rightarrow \omega_\pi$ and the support of the cokernel, $Z\subset
\mc{C}$, is a closed substack that is smooth and such that
$\pi:Z\rightarrow \mf{M}_{0,0}$ is unramified and is the normalization
of $\Delta$.  Over $U_1$, the lemma immediately follows from this and
the arguments in the proof of Proposition~\ref{prop-Qomega}.  As in
that case, it suffices to establish the lemma over $U_1$.
\end{proof}

\subsection{Computation of $Q_\pi(L)$ for invertible sheaves of degree
  $0$} \label{subsec-QL}
\marpar{subsec-QL}

\noindent
Let $M$ be an Artin stack, let $\pi:C\rightarrow M$ be a flat
1-morphism, relatively representable by proper algebraic spaces whose
geometric fibers are connected, at-worst-nodal curves of arithmetic
genus $0$.  Let $L$ be an invertible sheaf on $C$ of relative degree
$0$ over $M$.  This determines a morphism to the relative Picard of
the universal curve over $\mf{M}_{0,0}$, i.e., $\zeta_M:M \rightarrow
\mf{M}_{\ZZ}(\tau_{0,0}(0))$ such that the pullback of $\mc{C}$ is
equivalent to $C$, and such that the pullback of
$\OO_{\mc{C}}(\mc{D})$ differs from $L$ by $\pi^*L'$ for an invertible
sheaf $L'$ on $M$.  By Lemma~\ref{lem-pullback} and
Lemma~\ref{lem-inv}, $Q_\pi(L) \cong \zeta_M^*
Q_\pi(\OO_{\mc{C}}(\mc{D}))$.

\medskip\noindent
Let $\pi:\mc{C} \rightarrow \mf{M}_{\ZZ}(\tau_{0,0}(0))$ be the
universal curve.

\begin{prop} \label{prop-lcomp}
\marpar{prop-lcomp}
Over $\mf{M}_{\ZZ}(\tau_{0,0}(0))$, $\pi_*E_\pi(\mc{D}) = (0)$ and 
$R^1\pi_*
E_\pi(\mc{D})$ is a sheaf supported on $\Delta$.  The stalk of
$R^1\pi_* E_\pi(\mc{D})$ at the generic point of $\Delta_a$ is a
torsion sheaf of length $a^2$.  The filtration by order of vanishing
at the generic point has associated graded pieces of length $2a-1,
2a-3, \dots, 3,1$.  
\end{prop}

\begin{proof}
Over the open complement of $\Delta$, the divisor $\mc{D}$ is $0$.  So
the first part of the proposition reduces to the statement that
$R\pi_* E_\pi$ 
is quasi-isomorphic to $0$.  By definition of $E_\pi$, there is an
exact triangle,
$$
\begin{CD}
R\pi_* E_\pi @>>> R\pi_* \OO_{\mc{C}} @> \delta >> R\pi_* \omega_\pi
[1] @>>> 
R\pi_* E_\pi[1].
\end{CD}
$$
Of course the canonical isomorphism $R\pi_* \OO_{\mc{C}} \cong
\OO_{\mf{M}}$, and 
$E_\pi$ were defined so that the composition of $\delta$ with the trace
map, which is a quasi-isomorphism in this case, would be the
identity.  Therefore $\delta$ is a quasi-isomorphism, so $R\pi_*
E_\pi$ is quasi-isomorphic to $0$.  

\medskip\noindent
The second part can be proved, and to an extent only makes sense,
after smooth base-change to a scheme.  Let $\PP^1_s$ be a copy of
$\PP^1$ with homogeneous coordinates $S_0,S_1$.  Let $\PP^1_x$ be a
copy of $\PP^1$ with homogeneous coordinates $X_0,X_1$.  Let $\PP^1_y$
be a copy of $\PP^1$ with homogeneous coordinates $Y_0,Y_1$.  Denote
by $C \subset \PP^1_s \times \PP^1_x \times \PP^1_y$ the divisor with
defining equation $F=S_0X_0Y_0 - S_1X_1Y_1$.  The projection
$\text{pr}_s:C \rightarrow \PP^1_s$ is a proper, flat morphism whose
geometric fibers are connected, at-worst-nodal curves of arithmetic
genus $0$.  Denote by $L$ the invertible sheaf on $C$ that is the
restriction of $\text{pr}_x^*\OO_{\PP^1_x}(a)\otimes
\text{pr}_y^*\OO_{\PP^1_y}(-a)$.  This is an invertible sheaf of
relative degree $0$.  Therefore there is an induced $1$-morphism
$\zeta:\PP^1_s \rightarrow \mf{M}_{\ZZ}(\tau_{0,0}(0))$.  

\medskip\noindent
It is straightforward that $\zeta$ is smooth, and the image intersects
$\Delta_b$ iff $b=a$.  Moreover, $\zeta^* \Delta_a$ is the reduced
Cartier divisor $\mathbb{V}(S_0S_1) \subset \PP^1_s$.  There is an
obvious involution $i:\PP^1_s\rightarrow \PP^1_s$ by
$i(S_0,S_1)=(S_1,S_0)$, and $\zeta\circ i$ is $2$-equivalent to
$\zeta$.  Therefore the length of the $R^1\text{pr}_{s,*}
E_{\text{pr}_s}\otimes L$ is $2$ times the length of the stalk of
$R^1\pi_* E_{\pi}(\mc{D})$ at the generic point of $\Delta_a$; more
precisely, the length of the stalk at each of $(1,0),(0,1)\in \PP^1_s$
is the length of the stalk at $\Delta_a$.  Similarly for the lengths
of the associated graded pieces of the filtration.

\medskip\noindent
Because $E_{\text{pr}_s}$ is the extension class of the Trace mapping,
$R^1\text{pr}_{s,*}E_{\text{pr}_s} \otimes L$ is the cokernel of the
$\OO_{\PP^1_s}$-homomorphisms,
$$
\gamma:\text{pr}_{s,*}(L) \rightarrow
\text{Hom}_{\OO_{\PP^1_s}}(\text{pr}_{s,*}(L^\vee), \OO_{\PP^1_s}),
$$
induced via adjointness from the multiplication map,
$$
\text{pr}_{s,*}(L) \otimes \text{pr}_{s,*}(L^\vee) \rightarrow
\text{pr}_{s,*}(\OO_C) = \OO_{\PP^1_s}.
$$

\medskip\noindent
On $\PP^1_s\times \PP^1_x \times \PP^1_y$ there is a locally free
resolution of the push-forward of $L$, resp. $L^\vee$,
$$
\begin{array}{c}
0 \rightarrow \OO_{\PP^1_s}(-1)\boxtimes \OO_{\PP^1_x}(a-1)\boxtimes
\OO_{\PP^1_y}(-a-1) \xrightarrow{F} \OO_{\PP^1_s}(0)\boxtimes
\OO_{\PP^1_x}(a) \boxtimes \OO_{\PP^1_y}(-a) \rightarrow L \rightarrow
0, \\
0 \rightarrow \OO_{\PP^1_s}(-1)\boxtimes \OO_{\PP^1_x}(-a-1)\boxtimes
\OO_{\PP^1_y}(a-1) \xrightarrow{F} \OO_{\PP^1_s}(0)\boxtimes
\OO_{\PP^1_x}(-a) \boxtimes \OO_{\PP^1_y}(a) \rightarrow L^\vee
\rightarrow 0
\end{array}
$$
Hence $R\text{pr}_{s,*}L$ is the complex,
$$
\OO_{\PP^1_s}(-1)\otimes_k H^0(\PP^1_x,\OO_{\PP^1_x}(a-1)) \otimes_k
H^1(\PP^1_y,\OO_{\PP^1_y}(-a-1)) \xrightarrow{F} \OO_{\PP^1_s}
\otimes_k H^0(\PP^1_x,\OO_{\PP^1_x}(a))\otimes_k
H^1(\PP^1_y,\OO_{\PP^1_y}(-a)).
$$
Similarly for $R\text{pr}_{s,*}L^\vee$.  It is possible to write out
this map explicitly in terms of bases for $H^0$ and $H^1$, but for the
main statement just observe the complex has rank $1$ and degree $-a^2$.
Similarly for $R\text{pr}_{s,*}L^\vee$.  Therefore $R^1\pi_* E_\pi(L)$
is a torsion sheaf of length $2a^2$.  Because it is equivariant for
$i$, the localization at each of $(0,1)$ and $(1,0)$ has length $a^2$.

\medskip\noindent
The lengths of the associated graded pieces of the filtration by order
of vanishing at $\mathbb{V}(S_0S_1)$ can be computed from the
complexes for $R\text{pr}_{s,*}L$ and $R\text{pr}_{s,*}L^\vee$.  This
is left to the reader.
\end{proof}

\begin{cor} \label{cor-lcomp}
\marpar{cor-lcomp}
In the universal case, $Q_\pi(\mc{D}) = -\sum_{a\geq 0} a^2 \Delta_a$.
Therefore in the general case of $\pi:C\rightarrow M$ and an
invertible sheaf $L$ of relative degree $0$,
$$
Q_\pi(L) = {\sum_{\beta',\beta''}}^\prime \langle C_1(L),\beta'
\rangle \langle 
C_1(L), \beta'' \rangle \Delta_{\beta',\beta''}.
$$
\end{cor}

\section{Some divisor class relations} \label{sec-div}
\marpar{sec-div}

\noindent
In this section, Proposition~\ref{prop-Qomega} and
Proposition~\ref{prop-lcomp} are used to deduce several other divisor
class relations.  As usual, let $M$ be an Artin stack and let
$\pi:C\rightarrow M$ be a flat 1-morphism, relatively representable by
proper algebraic spaces whose geometric fibers are connected,
at-worst-nodal curves of genus $0$.

\begin{hyp} \label{hyp-GRR}
\marpar{hyp-GRR}
There are cycle class groups for $C$ and $M$ admitting Chern classes
for locally free sheaves, and such that Grothendieck-Riemann-Roch
holds for $\pi$.
\end{hyp}

\begin{lem} \label{lem-rel1}
For every Cartier divisor class $D$ on $C$ of relative degree $\langle
D, \beta \rangle$ over $M$, modulo $2$-power torsion,
$$
\pi_*(D\cdot D) + \langle D,\beta \rangle \pi_*(D\cdot
C_1(\omega_\pi)) = {\sum_{\beta',\beta''}}^\prime \langle D,\beta'
\rangle \langle D,\beta'' \rangle \Delta_{\beta',\beta''}.
$$
\end{lem}

\begin{proof}
Define $D' = 2D + \langle D,\beta \rangle C_1(\omega_\pi)$.  This is a
Cartier divisor class of relative degree $0$.  By
Corollary~\ref{cor-lcomp}, 
$$
Q_\pi(D') = {\sum_{\beta',\beta''}}^\prime (\langle 2D,\beta' \rangle
-\langle D,\beta \rangle)(\langle 2D,\beta'' \rangle - \langle D,\beta
\rangle) \Delta_{\beta',\beta''}.
$$
By Lemma~\ref{lem-GRR} this is,
$$
\begin{array}{c}
4\pi_*(D\cdot D) +4\langle D,\beta \rangle \pi_*(D\cdot
C_1(\omega_\pi) + (\langle D,\beta \rangle)^2 Q_\pi(C_1(\omega_\pi)) =
\\
{\sum_{\beta',\beta''}}^\prime (4\langle D,\beta' \rangle \langle
D,\beta'' \rangle - (\langle D,\beta \rangle)^2)
\Delta_{\beta',\beta''}.
\end{array}
$$
By Proposition~\ref{prop-Qomega}, $Q_\pi(\omega_\pi) =
-{\sum_{\beta',\beta''}}' \Delta_{\beta',\beta''}$.  Substituting this
into the equation, simplifying, and dividing by 4 gives the relation.
\end{proof}

\begin{lem} \label{lem-rel2}
\marpar{lem-rel2}
For every pair of Cartier divisor classes on $C$, $D_1,D_2$, of
relative degrees $\langle D_1,\beta \rangle$, resp. $\langle D_2,\beta
\rangle$, modulo $2$-power torsion,
$$
\begin{array}{c}
2\pi_*(D_1\cdot D_2) + \langle D_1,\beta \rangle \pi_*(D_2\cdot
C_1(\omega_\pi)) + \langle D_2,\beta \rangle \pi_*(D_1\cdot
C_1(\omega_\pi)) = 
\\
\\
{\sum_{\beta',\beta''}}^\prime (\langle D_1,\beta' \rangle \langle
D_2,\beta'' 
\rangle + \langle D_2,\beta' \rangle \langle D_1,\beta'' \rangle)
\Delta_{\beta',\beta''}.
\end{array}
$$
\end{lem}

\begin{proof}
This follows from Lemma~\ref{lem-rel1} and the polarization identity
for quadratic forms.
\end{proof}

\begin{lem} \label{lem-rel3}
\marpar{lem-rel3}
For every section of $\pi$, $s:M\rightarrow C$, whose image is
contained in the smooth locus of $\pi$,
$$
s(M)\cdot s(M) + s(M)\cdot C_1(\omega_\pi).
$$
\end{lem}

\begin{proof}
This follows by adjunction since the relative dualizing sheaf of
$s(M)\rightarrow M$ is trivial.
\end{proof}

\begin{lem} \label{lem-rel4}
\marpar{lem-rel4}
For every section of $\pi$, $s:M\rightarrow C$, whose image is
contained in the smooth locus of $\pi$ and for every Cartier divisor
class $D$ on $C$ of relative degree $\langle D,\beta \rangle$ over
$M$, modulo $2$-power torsion,
$$
\begin{array}{c}
2\langle D,\beta \rangle s^*D  -\pi_*(D\cdot D) - \langle D,\beta
\rangle^2 \pi_*(s(M)\cdot s(M)) = \\
\\
{\sum_{\beta',\beta''}}^\prime
(\langle D,\beta' \rangle^2 \langle s(M),\beta'' \rangle + \langle
D,\beta'' \rangle^2 \langle s(M),\beta' \rangle )
\Delta_{\beta',\beta''}.
\end{array}
$$
\end{lem}

\begin{proof}
By Lemma~\ref{lem-rel2},
$$
\begin{array}{c}
2s^* D + \pi_*(D\cdot C_1(\omega_\pi)) + \langle D,\beta \rangle
\pi_*(s(M)\cdot C_1(\omega_\pi)) = \\
\\
{\sum}^\prime (\langle D,\beta' \rangle \langle s(M),\beta'' \rangle +
\langle D,\beta'' \rangle \langle s(M),\beta' \rangle
)\Delta_{\beta',\beta''}.
\end{array}
$$
Multiplying both sides by $\langle D,\beta \rangle$,
$$
\begin{array}{c}
2\langle D,\beta \rangle s^* D + \langle D,\beta \rangle \pi_*(D \cdot
C_1(\omega_\pi)) + \langle D, \beta \rangle^2 \pi_*(s(M)\cdot
C_1(\omega_\pi)) = \\
\\
{\sum}^\prime(\langle D,\beta \rangle \langle D,\beta' \rangle \langle
s(M),\beta'' \rangle + \langle D,\beta \rangle \langle D,\beta''
\rangle \langle s(M),\beta' \rangle ) \Delta_{\beta',\beta''}.
\end{array}
$$
First of all, by Lemma~\ref{lem-rel4}, $\langle D,\beta \rangle^2
\pi_*(s(M) \cdot C_1(\omega_\pi)) = - \langle D,\beta \rangle^2
\pi_*(s(M)\cdot s(M))$.  Next, by Lemma~\ref{lem-rel1},
$$
\langle D,\beta \rangle \pi_*(D\cdot C_1(\omega_\pi)) = -\pi_*(D\cdot
D) + {\sum}^\prime \langle D,\beta' \rangle \langle D,\beta''
\Delta_{\beta',\beta''}.
$$
Finally,
$$
\begin{array}{c}
\langle D,\beta \rangle \langle D,\beta' \rangle \langle s(M),\beta''
\rangle + \langle D,\beta \rangle \langle D,\beta'' \rangle \langle
s(M),\beta' \rangle = \\ \\
(\langle D,\beta' \rangle + \langle D,\beta'' \rangle) \langle
D,\beta' \rangle \langle s(M),\beta'' \rangle + (\langle D,\beta'
\rangle + \langle D,\beta'' \rangle) \langle D,\beta'' \rangle \langle
s(M), \beta' \rangle = \\ \\
\langle D,\beta' \rangle^2 \langle s(M),\beta'' \rangle + \langle
D,\beta'' \rangle^2 \langle s(M),\beta' \rangle + \langle D,\beta'
\rangle \langle D,\beta'' \rangle (\langle s(M),\beta' \rangle +
\langle s(M),\beta'' \rangle) = \\ \\
\langle D,\beta' \rangle^2 \langle s(M),\beta'' \rangle + \langle
D,\beta'' \rangle^2 \langle s(M),\beta' \rangle + \langle D,\beta'
\rangle \langle D,\beta'' \rangle.
\end{array}
$$
Plugging in these 3 identities and simplifying gives the relation.
\end{proof}

\medskip\noindent
Let $\mc{C}$ be the universal curve over $\mf{M}_{0,0}$.  Let
$\mc{C}_\text{smooth}$ denote the smooth locus of $\pi$.  The
2-fibered product $\text{pr}_1:\mc{C}_\text{smooth}
\times_{\mf{M}_{0,0}} \mc{C} 
\rightarrow \mc{C}_\text{smooth}$ together with the diagonal
$\Delta:\mc{C}_\text{smooth} \rightarrow \mc{C}_\text{smooth}
\times_{\mf{M}_{0,0}} \mc{C}$ determine a 1-morphism
$\mc{C}_\text{smooth} \rightarrow \mf{M}_{1,0}$.  This extends to a
1-morphism $\mc{C} \rightarrow \mf{M}_{1,0}$.  The pullback of the
universal curve is a 1-morphism $\pi':\mc{C}' \rightarrow \mc{C}$ that
factors through $\text{pr}_1:\mc{C}\times_{\mf{M}_{0,0}} \mc{C}
\rightarrow \mc{C}$.  Denote the pullback of the universal section by
$s:\mc{C} \rightarrow \mc{C}'$.  Now $\mc{C}$ is regular, and the
complement of $\mc{C}_\text{smooth}$ has codimension $2$.  In
particular, $s^*\OO_{\mc{C}'}(s(\mc{C}))$ can be computed on
$\mc{C}_\text{smooth}$.  But the restriction to $\mc{C}_\text{smooth}$
is clearly $\omega^\vee_\pi$.  Therefore $s^*\OO_{\mc{C}'}(s(\mc{C}))
\cong 
\omega_\pi^\vee$ on all of $\mc{C}$.

\medskip\noindent
Pulling this back by $\zeta_C:C\rightarrow \mf{C}$ gives a 1-morphism
$\pi':C'\rightarrow C$ that factors through $\text{pr}_1:C\times_M
C\rightarrow C$.  Let $D$ be a Cartier divisor class on $C$ and
consider the pullback to $C'$ of $\text{pr}_2^* D$ on $C\times_M C$.
This is a Cartier divisor class $D'$ on $C'$.  Of course $s^* D' =
D$.  Moreover, by the projection formula the pushforward to $C\times_M
C$ of $D'\cdot D'$ is $\text{pr}_2^*(D\cdot D)$.  Therefore
$(\pi')_*(D'\cdot D')$ is $(\text{pr}_1)_*\text{pr}_2^*(D\cdot D)$,
i.e., $\pi^*\pi_*(D\cdot D)$.  Finally, denote by,
$$
\sum_{\beta',\beta''} \langle D,\beta'' \rangle^2
\widetilde{\Delta}_{\beta',\beta''},
$$
the divisor class on $C$,
$$
{\sum_{\beta',\beta''}}^\prime (\langle D,\beta'' \rangle^2 \langle
s,\beta' \rangle + \langle D,\beta' \rangle^2 \langle s,\beta''
\rangle)\Delta_{\beta',\beta''}.
$$
The point is this: if $\pi$ is smooth over every generic point of $M$,
then the divisor class $\widetilde{\Delta}_{\beta',\beta''}$ is the
irreducible component of $\pi^{-1}(\Delta_{\beta',\beta''})$
corresponding to the vertex $v'$, i.e., the irreducible component with
``curve class'' $\beta'$.
Putting this all together and applying Lemma~\ref{lem-rel4} gives the
following.

\begin{lem} \label{lem-rel5}
\marpar{lem-rel5}
For every Cartier divisor class $D$ on $C$ of relative degree $\langle
D,\beta \rangle$ over $M$,
$$
\begin{array}{c}
2\langle D,\beta \rangle D - \pi^*\pi_*(D\cdot D) + \langle D,\beta
\rangle^2 C_1(\omega_\pi) = \\ \\
\sum_{\beta',\beta''} \langle D,\beta'' \rangle^2
\widetilde{\Delta}_{\beta',\beta''}.
\end{array}
$$
In particular, the relative Picard group of $\pi$ is generated by
$C_1(\omega_\pi)$ and the boundary divisor classes
$\widetilde{\Delta}_{\beta',\beta''}$.  
\end{lem}

\begin{rmk} \label{rmk-rel5}
\marpar{rmk-rel5}
If $\langle D,\beta \rangle \neq 0$ then, at least up to torsion,
Lemma~\ref{lem-rel1} follows from Lemma~\ref{lem-rel5} by intersecting
both sides of the relation by $D$ and then applying $\pi_*$.  This was
pointed out by Pandharipande, who also proved Lemma~\ref{lem-rel4} up
to numerical equivalence in ~\cite[Lem. 2.2.2]{QDiv} (by a very
different method).   
\end{rmk}

\begin{lem} \label{lem-rel6}
\marpar{lem-rel6}
Let $s,s':M\rightarrow C$ be sections with image in the smooth locus
of $\pi$ such that $s(M)$ and $s'(M)$ are disjoint.  Then,
$$
\pi_*(s(M)\cdot s(M)) + \pi_*(s'(M)\cdot s'(M)) = -
\sum_{\beta',\beta''} \langle s(M),\beta' \rangle \langle
s'(M),\beta'' \rangle \Delta_{\beta',\beta''}.
$$
\end{lem}

\begin{proof}
Apply Lemma~\ref{lem-rel2} and use $s(M)\cdot s'(M)=0$ and
Lemma~\ref{lem-rel3}.  
\end{proof}

\begin{lem} \label{lem-rel7}
Let $r\geq 2$ and $s_1,\dots,s_r:M\rightarrow C$ be sections with image in the
smooth locus of $\pi$ and which are pairwise disjoint.  Then,
$$
-\sum_{i=1}^r \pi_*(s_i(M)\cdot s_i(M)) = (r-2)\pi_*(s_1(M)\cdot
s_1(M)) + \sum_{\beta',\beta''} \langle s_1(M),\beta' \rangle \langle
s_2(M) + \dots + s_r(M),\beta'' \rangle \Delta_{\beta',\beta''}.
$$
\end{lem}

\begin{proof}
This follows from Lemma~\ref{lem-rel6} by induction.
\end{proof}

\begin{lem} \label{lem-rel8}
\marpar{lem-rel8}
Let $r\geq 2$ and let $s_1,\dots,s_r:M \rightarrow C$ be sections with
image in the smooth locus of $\pi$ and which are pairwise disjoint.
Then,
$$
-\sum_{i=1}^r \pi_*(s_i(M)\cdot s_i(M)) = r(r-2) \pi_*(s_1(M)\cdot
s_1(M)) + \sum_{\beta',\beta''} \langle s_1(M),\beta' \rangle \langle
s_2(M) + \dots + s_r(M),\beta'' \rangle^2 \Delta_{\beta',\beta''}.
$$
Combined with Lemma~\ref{lem-rel7} this gives,
$$
\begin{array}{c}
(r-1)(r-2)\pi_*(s_1(M)\cdot s_1(M)) = \\ \\
-\sum_{\beta',\beta''} \langle
s_1(M),\beta' \rangle \langle s_2(M) + \dots + s_r(M),\beta'' \rangle
(\langle s_2(M)+\dots +s_r(M), \beta'' \rangle - 1)
\Delta_{\beta',\beta''},
\end{array}
$$
which in turn gives,
$$
\begin{array}{c}
-(r-1)\sum_{i=1}^r\pi_*(s_i(M)\cdot s_i(M)) = \\ \\
\sum_{\beta',\beta''} \langle s_1(M),\beta' \rangle \langle s_2(M) +
\dots + s_r(M), \beta'' \rangle (r- \langle s_2(M) +\dots + s_r(M),
\beta'' \rangle) \Delta_{\beta',\beta''}.
\end{array}
$$
In the notation of Example~\ref{ex-boundary}, this is,
$$
-(r-1)(r-2)\pi_*(s_1(M)\cdot s_1(M)) = \sum_{(A,B), \ 1\in A}
\#B(\#B-1) \Delta_{(A,B)},
$$
and
$$
-(r-1)\sum_{i=1}^r \pi_*(s_i(M)\cdot s_i(M)) = \sum_{(A,B), \ 1\in A}
\#B(r-\#B) \Delta_{(A,B)}.
$$
\end{lem}

\begin{proof}
Denote $D=\sum_{i=2}^r s_i(M)$.  Apply Lemma~\ref{lem-rel4} to get,
$$
\begin{array}{c}
2(r-1)\cdot 0 - \sum_{i=2}^r \pi_*(s_i(M)\cdot s_i(M)) -
(r-1)^2\pi_*(s_1(M)\cdot s_1(M)) = \\ \\
\sum_{\beta',\beta''} \langle s_1(M),\beta' \rangle \langle
s_2(M)+\dots + s_r(M),\beta'' \rangle^2 \Delta_{\beta',\beta''}.
\end{array}
$$
Simplifying,
$$
-\sum_{j=1}^r \pi_*(s_i(M)\cdot s_i(M)) = r(r-2) \pi_*(s_1(M)\cdot
s_1(M)) + \sum \langle s_1(M),\beta' \rangle \langle s_2(M)+ \dots
+s_r(M),\beta'' \rangle^2 \Delta_{\beta',\beta''}.
$$
Subtracting from the relation in Lemma~\ref{lem-rel7} gives the
relation for $(r-1)(r-2)\pi_*(s_1(M)\cdot s_1(M))$.  
Multiplying the first relation by $(r-1)$, plugging in the second
relation and simplifying gives the third relation.
\end{proof}

\begin{lem} \label{lem-rel9}
\marpar{lem-rel9}
Let $r\geq 2$ and let $s_1,\dots,s_r:M\rightarrow C$ be everywhere
disjoint sections with image in the smooth locus.  For every $1\leq i
< j \leq r$, using the notation from Example~\ref{ex-boundary},
$$
\sum_{(A,B), \ i\in A} \#B(r-\#B) \Delta_{(A,B)} = \sum_{(A',B'), j \in
  A} \#B'(r-\#B') \Delta_{(A',B')}.
$$
\end{lem}

\begin{proof}
This follows from Lemma~\ref{lem-rel8} by permuting the roles of $1$
with $i$ and $j$.
\end{proof}

\begin{lem} \label{lem-rel10}
\marpar{lem-rel10}
Let $r \geq 2$ and let $s_1,\dots,s_r:M \rightarrow C$ be everywhere
disjoint sections with image in the smooth locus of $\pi$.  For every
Cartier 
divisor class $D$ on $C$ of relative degree $\langle D,\beta \rangle$,
$$
2(r-1)(r-2)\langle D,\beta \rangle s_1^*D = (r-1)(r-2)\pi_*(D\cdot D)
+ \sum_{\beta',\beta''} \langle s_1(M),\beta' \rangle a(D,\beta'')
\Delta_{\beta',\beta''},
$$
where, 
$$
a(D,\beta'') = (r-1)(r-2)\langle D,\beta'' \rangle^2 -
\langle D,\beta \rangle^2 \langle s_2(M)+\dots +s_r(M), \beta'' \rangle (
\langle s_2(M)+\dots + s_r(M),\beta'' \rangle -1).
$$
In particular, if $r\geq 3$, then modulo torsion $s_i^*D$ is in the
span of $\pi_*(D\cdot D)$ and boundary divisors for every
$i=1,\dots,r$.
\end{lem}

\begin{proof}
This follows from Lemma~\ref{lem-rel4} and Lemma~\ref{lem-rel8}.
\end{proof}

\begin{lem} \label{lem-rel11}
\marpar{lem-rel11}
Let $r\geq 2$ and let $s_1,\dots,s_r:M\rightarrow C$ be everywhere
disjoint sections with image in the smooth locus of $\pi$.  Consider
the sheaf $\mc{E} = \Omega_\pi(s_1(M)+\dots+s_r(M))$.  The perfect
complex $R\pi_*R\textit{Hom}_{\OO_C}(\mc{E},\OO_C)$ has rank $3-r$ and
the first Chern class of the determinant is $-2\Delta
-\sum_{i=1}^r(s_i(M)\cdot s_i(M))$.  In particular, if $r\geq 2$, up
to torsion,
$$
\begin{array}{c}
C_1(\text{det}R\pi_*
R\textit{Hom}_{\OO_C}(\Omega_\pi(s_1(M)+\dots+s_r(M)),\OO_C)) = \\ \\
-2\Delta + \frac{1}{r-1}\sum_{(A,B),\ 1\in A} \#B(r-\#B)
\Delta_{(A,B)}.
\end{array}
$$
\end{lem}

\begin{proof}
There is a short exact sequence,
$$
\begin{CD}
0 @>>> \Omega_\pi @>>> \Omega_\pi(s_1(M)+\dots+s_r(M))
@>>> \oplus_{i=1}^r (s_i)_*\OO_M @>>> 0.
\end{CD}
$$
Combining this with Lemma~\ref{lem-Omega}, Lemma~\ref{lem-rel8}, and
chasing through exact sequences gives the lemma.
\end{proof}

\section{The virtual canonical bundle} \label{sec-vircan}
\marpar{sec-vircan}

\noindent
Let $k$ be a field, let $X$ be a connected, smooth algebraic space
over $k$ of dimension $n$, let
$M$ be an Artin stack over $k$, let $\pi:C\rightarrow M$ be a flat
1-morphism, 
relatively representable by proper algebraic spaces whose geometric
fibers are connected, at-worst-nodal curves of arithmetic genus $0$,
let $s_1,\dots,s_r:M\rightarrow C$ be pairwise disjoint sections with
image contained in the smooth locus of $\pi$ (possibly $r=0$, i.e.,
there are no sections), and let $f:C\rightarrow
X$ be a 1-morphism of $k$-stacks.  In this setting, Behrend and
Fantechi introduced a perfect complex $E^\bullet$ on $M$ of amplitude
$[-1,1]$ and a morphism to the cotangent complex,
$\phi:E^\bullet \rightarrow L_M^\bullet$, ~\cite{BM}.  
If $\text{char}(k)=0$ and $M$ is the Deligne-Mumford stack of stable
maps to $X$, Behrend and Fantechi prove $E^\bullet$ has amplitude
$[-1,0]$,
$h^0(\phi)$ is an isomorphism
and $h^{-1}(\phi)$ is surjective.  In
many interesting cases, $\phi$ is a quasi-isomorphism.  Then
$\text{det}(E^\bullet)$ is an invertible 
dualizing sheaf for $M$. 
Because of this, 
$\text{det}(E^\bullet)$ is called the \emph{virtual canonical
  bundle}.  In this section the relations from Section~\ref{sec-div}
are used to give a formula for the divisor class of the virtual
canonical bundle.  Hypothesis ~\ref{hyp-GRR} holds for $\pi$.

\medskip\noindent
Denote by $L_{(\pi,f)}$ the cotangent complex of the morphism
$(\pi,f):C\rightarrow M\times X$.  This is a perfect complex of
amplitude $[-1,0]$.  There is a distinguished triangle,
$$
\begin{CD}
L_\pi @>>> L_{(\pi,f)} @>>> f^* \Omega_X[1] @>>> L_\pi[1].
\end{CD}
$$
There is a slight variation $L_{(\pi,f,s)}$ taking into account the
sections which fits into a distinguished triangle,
$$
\begin{CD}
L_\pi(s_1(M)+\dots+s_r(M)) @>>> L_{(\pi,f,s)} @>>> f^*\Omega_X[1] @>>>
L_\pi(s_1(M)+\dots+s_r(M))[1]. 
\end{CD}
$$
The complex $E^\bullet$ is defined
to be $(R\pi_*(L_{(\pi,f,s)}^\vee)[1])^\vee$, where $(F^\bullet)^\vee$
is $R\textit{Hom}(F^\bullet,\OO)$.  In particular,
$\text{det}(E^\bullet)$ is the determinant of
$R\pi_*(L_{(\pi,f,s)}^\vee)$.  From the distinguished triangle,
$\text{det}(E^\bullet)$ is
$$
\text{det}(R\pi_*
R\textit{Hom}_{\OO_C}(\Omega_\pi(s_1(M)+\dots+s_r(M)),\OO_C)) 
\otimes \text{det}(R\pi_* f^*T_X)^\vee.
$$
By Lemma~\ref{lem-rel11}, the first term is known.  
The second term
follows easily from Grothendieck-Riemann-Roch.

\begin{lem} \label{lem-TX}
\marpar{lem-TX}
Assume that the relative degree of $f^*C_1(\Omega_X)$ is nonzero.
Then $R\pi_*f^*T_X[-1]$ has rank $\langle -f^*C_1(\Omega_X),\beta
\rangle + n$, and up to torsion the first Chern class of the
determinant is,
$$
\begin{array}{c}
\frac{1}{2\langle -f^*C_1(\Omega_X),\beta \rangle} \lt[ 2\langle
-f^*C_1(\Omega_X), \beta \rangle \pi_* f^* C_2(\Omega_X) \rt. \\
\\
 - (\langle
-f^* C_1(\Omega_X),\beta \rangle + 1) \pi_* f^* C_1(\Omega_X)^2 +
\\ \\
\lt.
{\sum}' \langle -f^*C_1(\Omega_X),\beta' \rangle \langle
-f^*C_1(\Omega_X),\beta'' \rangle \Delta_{\beta',\beta''}\rt].
\end{array}
$$
\end{lem}

\begin{proof}
The Todd class $\tau_\pi$ of $\pi$ is $1 - \frac{1}{2}C_1(\omega_\pi)
+ \tau_2 + \dots$, 
where $\pi_*\tau_2 = 0$.  The Chern character of $f^*T_X$ is,
$$
n - f^*C_1(\Omega_X) + \frac{1}{2}(f^*C_1(\Omega_X)^2 -
2f^*C_2(\Omega_X)) + \dots
$$
Therefore $\text{ch}(f^*T_X)\cdot \tau_\pi$ equals,
$$
n -
\lt[f^*C_1(\Omega_X)+\frac{n}{2}C_1(\Omega_\pi) \rt] + \frac{1}{2}\lt[
f^*C_1(\Omega_X)^2 - 2f^*C_2(\Omega_X) + f^*C_1(\Omega_X)\cdot
C_1(\omega_\pi) \rt] + n\tau_2 + \dots
$$
Applying $\pi_*$ and using that $\pi_*\tau_2 = 0$, the rank is
$n+\langle -f^*C_1(\Omega_X),\beta \rangle$, and the determinant has
first Chern class,
$$
\frac{1}{2}\pi_*\lt[ f^*C_1(\Omega_X)^2 - 2f^*C_2(\Omega_X) \rt] +
\frac{1}{2} \pi_*(f^*C_1(\Omega_X)\cdot C_1(\omega_\pi)).
$$
Applying Lemma~\ref{lem-rel1} and simplifying gives the relation.
\end{proof}  

\begin{prop} \label{prop-vc}
\marpar{prop-vc}
The rank of $E^\bullet$ is $\langle -f^*C_1(\Omega_X),\beta \rangle +
n + r - 3.$  The following divisor class reations hold modulo torsion.
If $\langle -f^*C_1(\Omega_X),\beta \rangle\neq 0$ and $r=0$,
the first Chern class of the virtual canonical bundle
is,
\begin{equation}
\begin{array}{c}
\frac{1}{2\langle -f^*C_1(\Omega_X),\beta \rangle} \lt[ 2\langle
-f^*C_1(\Omega_X), \beta \rangle \pi_* f^* C_2(\Omega_X) \rt. \\
\\
- (\langle
-f^* C_1(\Omega_X),\beta \rangle + 1) \pi_* f^* C_1(\Omega_X)^2 +
\\ \\
\lt.
{\sum}' (\langle -f^*C_1(\Omega_X),\beta' \rangle \langle
-f^*C_1(\Omega_X),\beta'' \rangle -4\langle -f^*C_1(\Omega_X),\beta
\rangle) \Delta_{\beta',\beta''}\rt].
\end{array}
\end{equation}
If $\langle -f^*C_1(\Omega_X),\beta \rangle \neq 0$ and $r=1$, 
the first Chern class of the virtual canonical bundle is,
\begin{equation}
\begin{array}{c}
\frac{1}{2\langle -f^*C_1(\Omega_X),\beta \rangle} \lt[ 2\langle
-f^*C_1(\Omega_X), \beta \rangle \pi_* f^* C_2(\Omega_X) \rt. \\
\\
- (\langle
-f^* C_1(\Omega_X),\beta \rangle + 1) \pi_* f^* C_1(\Omega_X)^2 +
\\ \\
\lt.
{\sum}' (\langle -f^*C_1(\Omega_X),\beta' \rangle \langle
-f^*C_1(\Omega_X),\beta'' \rangle -4\langle -f^*C_1(\Omega_X),\beta
\rangle) \Delta_{\beta',\beta''}\rt] \\
\\
-\pi_*(s_1(M)\cdot s_1(M)).
\end{array}
\end{equation}
If $\langle -f^*C_1(\Omega_X),\beta \rangle \neq 0$ and $r\geq 2$, 
the first Chern class of the virtual canonical bundle is,
\begin{equation}
\begin{array}{c}
\frac{1}{2\langle -f^*C_1(\Omega_X),\beta \rangle} \lt[ 2\langle
-f^*C_1(\Omega_X), \beta \rangle \pi_* f^* C_2(\Omega_X) \rt. \\
\\
- (\langle
-f^* C_1(\Omega_X),\beta \rangle + 1) \pi_* f^* C_1(\Omega_X)^2 +
\\ \\
\lt.
{\sum}' (\langle -f^*C_1(\Omega_X),\beta' \rangle \langle
-f^*C_1(\Omega_X),\beta'' \rangle -4\langle -f^*C_1(\Omega_X),\beta
\rangle) \Delta_{\beta',\beta''}\rt] \\
\\
+ \frac{1}{r-1} \sum_{(A,B), 1\in
  A} \#B(r-\#B)\Delta_{(A,B)}.
\end{array}
\end{equation}
\end{prop}

\bibliography{my}
\bibliographystyle{abbrv}

\end{document}